\documentclass[a4paper,11pt]{article}
\usepackage{latexsym}
\usepackage{amssymb}
\usepackage{enumerate}
\usepackage{amsfonts}
\usepackage{amsmath}
\usepackage[dvips]{graphicx}
\usepackage[paper=a4paper,left=30mm,right=20mm,top=25mm,bottom=30mm]{geometry}

\newcommand{\qed}{$\Box$}

\newenvironment{@abssec}[1]{%
    \if@twocolumn

      \section*{#1}%
    \else

      \vspace{.05in}\footnotesize
      \parindent .2in
 {\upshape\bfseries #1. }\ignorespaces
    \fi}

    {\if@twocolumn\else\par\vspace{.1in}\fi}

\newenvironment{keywords}{\begin{@abssec}{\keywordsname}}{\end{@abssec}}

\newenvironment{AMS}{\begin{@abssec}{\AMSname}}{\end{@abssec}}

\newcommand\keywordsname{Key words}
\newcommand\AMSname{AMS subject classifications}
\newcommand\AMname{AMS subject classification}
\newtheorem{theorem}{Theorem}
 \newtheorem{lemma}[theorem]{Lemma}
 \newtheorem{corollary}[theorem]{Corollary}
 
\newtheorem{remark}[theorem]{Remark}
 
\def\qed{\vbox{\hrule height0.6pt\hbox{%
  \vrule height1.3ex width0.6pt\hskip0.8ex
  \vrule width0.6pt}\hrule height0.6pt
 }}

\title{Two-phase heat conductors with a stationary isothermic surface and their related elliptic overdetermined problems
\thanks{This research was partially supported by the Grants-in-Aid
for Scientific Research (B) ($\sharp$ 26287020)  and Challenging Exploratory Research ($\sharp$ 16K13768) of
Japan Society for the Promotion of Science.}}

\author{Shigeru Sakaguchi\thanks{Research Center for Pure and Applied Mathematics,
Graduate School of  Information Sciences, Tohoku
University, Sendai, 980-8579,  Japan.
({\tt sigersak@tohoku.ac.jp}).}}

\begin{document}

\maketitle

\begin{abstract}
We consider a two-phase heat conductor in two dimensions consisting of a core and a shell with different constant conductivities. When the medium outside the two-phase conductor has a possibly different conductivity,  we consider the Cauchy problem in two dimensions where initially the conductor has temperature 0 and  the outside medium has temperature 1.
 It is shown that, if  there is a stationary isothermic surface in the shell near the boundary, then the structure of the conductor must be circular.  Moreover, as by-products of the method of the proof, we mention other proofs of all the previous results of \cite{Strieste2016} in $N(\ge 2)$ dimensions and two theorems on their related two-phase elliptic overdetermined problems.
  \end{abstract}

\begin{keywords}
heat equation, diffusion equation, two-phase heat conductor, transmission condition, Cauchy problem, two dimensions, stationary isothermic surface, symmetry, elliptic overdetermined problem
\end{keywords}

\begin{AMS}
Primary 35K05 ; Secondary  35K10,  35B06, 35B40,  35K15, 35K20, 35J05, 35J25
\end{AMS}

\pagestyle{plain}
\thispagestyle{plain}

\section{Introduction}
\label{introduction}

\vskip 2ex
Let $\Omega$ be a bounded $C^2$ domain in $\mathbb R^N\ (N \ge 2)$ with boundary $\partial\Omega$, and let $D$ be a bounded  $C^2$ open set in $\mathbb R^N$ which may have finitely many connected components.  Assume that $\Omega\setminus\overline{D}$ is connected and $\overline{D} \subset \Omega$. Denote by $\sigma=\sigma(x)\ (x \in \mathbb R^N)$  the conductivity distribution of the medium given by
$$
\sigma =
\begin{cases}
\sigma_c \quad&\mbox{in } D, \\
\sigma_s \quad&\mbox{in } \Omega \setminus D, \\
\sigma_m \quad &\mbox{in } \mathbb R^N \setminus \Omega,
\end{cases}
$$
where $\sigma_c, \sigma_s, \sigma_m$ are positive constants and $\sigma_c \not=\sigma_s$. This kind of three-phase electrical conductor has been dealt with in \cite{KLS2016} in the study of neutrally coated inclusions.

In the previous paper \cite{Strieste2016}, we considered the heat diffusion over two-phase or three-phase heat conductors.
Let $u =u(x,t)$ be the unique bounded solution of either the initial-boundary value problem for the diffusion equation:
\begin{eqnarray}
&&u_t =\mbox{ div}( \sigma \nabla u)\quad\mbox{ in }\ \Omega\times (0,+\infty), \label{heat equation initial-boundary}
\\
&&u=1  \ \quad\qquad\qquad\mbox{ on } \partial\Omega\times (0,+\infty), \label{heat Dirichlet}
\\ 
&&u=0  \  \quad\qquad\qquad \mbox{ on } \Omega\times \{0\},\label{heat initial}
\end{eqnarray}
or the Cauchy problem for the diffusion equation:
\begin{equation}
  u_t =\mbox{ div}(\sigma \nabla u)\quad\mbox{ in }\  \mathbb R^N\times (0,+\infty) \ \mbox{ and }\ u\ ={\mathcal X}_{\Omega^c}\ \mbox{ on } \mathbb R^N\times
\{0\},\label{heat Cauchy}
\end{equation}
where ${\mathcal X}_{\Omega^c}$ denotes the characteristic function of the set $\Omega^c=\mathbb R^N \setminus\Omega$. Consider a bounded domain $G$ in $\mathbb R^N$ satisfying
\begin{equation}
\label{near the boundary}
\overline{D} \subset G \subset \overline{G} \subset \Omega\ \mbox{ and } \mbox{ dist}(x,\partial\Omega) \le \mbox{ dist}(x, \overline{D})\ \mbox{ for every } x \in \partial G.
\end{equation}
In \cite{Strieste2016}, we obtained the following theorems.
\begin{theorem}[\cite{Strieste2016}]
\label{th:stationary isothermic} Let $u$ be the solution of problem \eqref{heat equation initial-boundary}-\eqref{heat initial} for $N \ge 2$, and let
$\Gamma$ be a connected component of $\partial G$ satisfying
\begin{equation}
\label{nearest component}
\mbox{\rm dist}(\Gamma, \partial\Omega) = \mbox{\rm dist}(\partial G, \partial\Omega).
\end{equation}  
If there exists a function $a : (0, +\infty) \to (0, +\infty) $ satisfying
\begin{equation}
\label{stationary isothermic surface partially}
u(x,t) = a(t)\ \mbox{ for every } (x,t) \in \Gamma \times (0, +\infty),
\end{equation}
then $\Omega$ and $D$ must be concentric balls.
\end{theorem}

\begin{corollary} [\cite{Strieste2016}]
\label{cor:stationary isothermic} Let $u$ be the solution of problem \eqref{heat equation initial-boundary}-\eqref{heat initial} for $N \ge 2$.  If
there exists a function $a : (0, +\infty) \to (0, +\infty) $ satisfying
\begin{equation}
\label{stationary isothermic surface}
u(x,t) = a(t)\ \mbox{ for every } (x,t) \in \partial G\times (0, +\infty),
\end{equation}
then $\Omega$ and $D$ must be concentric balls.
\end{corollary}

\begin{theorem}[\cite{Strieste2016}] 
\label{th:stationary isothermic cauchy} Let $u$ be the solution of problem \eqref{heat Cauchy} for $N \ge 3$. Then the following assertions hold:
\begin{itemize}
\item[\rm (a)] If
there exists a function $a : (0, +\infty) \to (0, +\infty) $ satisfying \eqref{stationary isothermic surface},
then $\Omega$ and $D$ must be concentric balls.
\item[\rm (b)] If $\sigma_s=\sigma_m$ and there exists a function $a : (0, +\infty) \to (0, +\infty) $ satisfying \eqref{stationary isothermic surface partially} for a connected component $\Gamma$ of $\partial G$ with \eqref{nearest component}, then  $\Omega$ and $D$ must be concentric balls.
\end{itemize}
\end{theorem}

In \cite{Strieste2016}, Theorem \ref{th:stationary isothermic cauchy} is limited to the case where $N \ge 3$, which is not natural; that is required for  technical reasons in the use of the auxiliary functions $U, V, W$ given in \cite[Proof of Theorem 1.3, pp. 184--186]{Strieste2016}. We conjectured that Theorem \ref{th:stationary isothermic cauchy} holds true also for $N=2$. 

The main purpose of the present paper is to show that this conjecture is true. Namely, we show the following theorem.

\begin{theorem}
\label{th:stationary isothermic cauchy in two dimensions} Let $u$ be the solution of problem \eqref{heat Cauchy} in two dimensions. Then the following assertions hold:
\begin{itemize}
\item[\rm (a)] If
there exists a function $a : (0, +\infty) \to (0, +\infty) $ satisfying \eqref{stationary isothermic surface},
then $\Omega$ and $D$ must be concentric disks.
\item[\rm (b)] If $\sigma_s=\sigma_m$ and there exists a function $a : (0, +\infty) \to (0, +\infty) $ satisfying \eqref{stationary isothermic surface partially} for a connected component $\Gamma$ of $\partial G$ with \eqref{nearest component}, then  $\Omega$ and $D$ must be concentric disks.
\end{itemize}
\end{theorem}

The other purpose  is to mention that the method employed in the present paper enables us to give other proofs of all the previous results of \cite{Strieste2016} in $N(\ge 2)$ dimensions and two theorems  on their related two-phase elliptic overdetermined problems (see section \ref{section5}). 

\vskip 2ex
The following sections are organized as follows. In section \ref{section2}, in two dimensions we give three preliminaries dealt with in \cite{Strieste2016} for the sake of convenience. 
Section \ref{section3} introduces four key tools concerning partial and ordinary differential equations. These four tools are stated in $N (\ge 2)$ dimensions.  In section \ref{section4}, we prove Theorem \ref{th:stationary isothermic cauchy in two dimensions}. If $D$ is not a disk, we use the transmission condition \eqref{transmission condition between U and V} on $\partial D$ to get a contradiction to either Hopf's boundary point lemma or Lemma \ref{le: the unique determination} stating the unique determination of the inclusions with one Cauchy data.  New auxiliary functions $U, V, W$  given in section \ref{section4}  play a key role. In section \ref{section5},  as by-products of the method of the proof in section \ref{section4}, we mention that we may have other proofs of all the previous results of \cite{Strieste2016} in $N(\ge 2)$ dimensions and two theorems (see Theorems \ref{th:constant Neumann boundary condition} and \ref{th:constant Dirichlet boundary condition on the inner sphere}) on related two-phase elliptic overdetermined problems.
 Indeed, the method of the proof employed in the present paper also gives other proofs of Theorems \ref{th:stationary isothermic} and \ref{th:stationary isothermic cauchy}. These new proofs do not use the explicit radially symmetric solutions of Poisson's equation over balls. See the radially symmetric solution $v =v(r)$ given in Remark \ref{poisson equation case} in section \ref{section3}. 
The proofs in $N (\ge 3)$ dimensions are 
parallel to that in two dimensions, since the four key tools are given in $N (\ge 2)$ dimensions and the preliminaries similar to those in section \ref{section2} are given in \cite[Section 2, pp. 169--180]{Strieste2016}. 
\setcounter{equation}{0}
\setcounter{theorem}{0}

\section{Preliminaries}
\label{section2}

Concerning the behavior of the solution of  problem \eqref{heat Cauchy} in two dimensions, we start with the following lemma.
\begin{lemma} 
\label{le:initial behavior exponential decay and decay at infinity} Let $u$ be the solution of problem \eqref{heat Cauchy} in two dimensions. We have the following assertions:
\begin{itemize}
\item[\rm (a)] 
$
0 < 1-u < 1 \quad \mbox{ in } \mathbb R^2 \times (0,\infty).
$
\item[\rm (b)] 
 $
 \lim\limits_{|x| \to \infty} (1-u(x,t)) = 0  \quad \mbox{ for every } t \in (0,\infty).
 $
 \end{itemize}
\end{lemma}

\noindent
{\it Proof.\ } 
We make use of the Gaussian bounds for the fundamental solutions of parabolic equations due to
Aronson\cite[Theorem 1, p. 891]{Ar1967bams}(see also \cite[p. 328]{FaS1986arma}). Let $g = g(x,\xi,t)$ be the fundamental solution of $u_t=\mbox{ div}(\sigma\nabla u)$. Then there exist two positive constants $\alpha$ and $M$ such that
\begin{equation}
\label{Gaussian bounds}
M^{-1}t^{-1}e^{-\frac{\alpha|x-\xi|^2}{t}}\le g(x,\xi,t) \le M t^{-1}e^{-\frac{|x-\xi|^2}{\alpha t}} 
\end{equation}
 for all $(x,t), (\xi,t) \in \mathbb R^2\times(0,+\infty)$. 
 For the solution $u$ of problem \eqref{heat Cauchy},  $1-u$ is regarded as the unique bounded solution of the Cauchy problem for the diffusion equation with initial data ${\mathcal X}_{\Omega}$. Hence we have from \eqref{Gaussian bounds}
 $$
 1-u(x,t) = \int_{\mathbb R^2}  g(x,\xi,t){\mathcal X}_{\Omega}(\xi)\ d\xi \le M t^{-1} \int_{\Omega} e^{-\frac{|x-\xi|^2}{\alpha t}}\ d\xi,
 $$
 which yields (b) and (a),  since  $g = g(x,\xi,t)$ is the fundamental solution.
 \qed
 
 Let us quote the following two lemmas from \cite[Lemma 2.4 and Lemma 2.5, p. 176 and p. 179]{Strieste2016} only for the Cauchy problem in two dimensions:

\begin{lemma} [\cite{Strieste2016}]
\label{le: constant weingarten curvature}
Let $u$ be the solution of problem \eqref{heat Cauchy} in two dimensions. 
Under the assumption \eqref{stationary isothermic surface partially}, the following assertions hold:
\begin{enumerate}[\rm (1)]
\item There exists a number $R > 0$ such that 
$$
\mbox{\rm dist}(x, \partial\Omega) = R\ \mbox{ for every } x \in \Gamma.
$$
\item $\Gamma$ is a real analytic regular curve.
\item There exists a connected component $\gamma$ of $\partial\Omega$, that is also a real analytic regular curve, such that the mapping $\gamma \ni y \mapsto x(y)  \equiv y-R\nu(y) \in \Gamma$, where $\nu(y)$ is the unit outward normal vector to $\partial\Omega$ at $y \in \gamma$,  is a diffeomorphism; in particular $\gamma$ and $\Gamma$  are parallel regular curves at distance $R$.
\item It holds that
\begin{equation}
\label{bounds of curvatures}
 \kappa(y) < \frac 1R\ \mbox{ for every } y \in \gamma,
\end{equation}
where $\kappa(y)$ is the curvature of $\partial\Omega$ at $y \in \gamma$ with respect to the unit inward normal vector $-\nu(y)$ to $\partial\Omega$.
\item There exists a number $c > 0$ such that
\begin{equation}
\label{monge-ampere}
\frac 1R-\kappa(y) = c\quad\mbox{ for every } y \in \gamma.
\end{equation}
\end{enumerate}
\end{lemma}

\begin{lemma}[\cite{Strieste2016}] 
\label{le: constant weingarten curvature 2}
Let $u$ be the solution of  problem \eqref{heat Cauchy} in two dimensions. 
Under the assumption \eqref{stationary isothermic surface}, the same assertions {\rm (1)--(5)} as in {\rm Lemma \ref{le: constant weingarten curvature}} hold provided $\Gamma$ and $\gamma$ are replaced by $\partial G$ and $\partial\Omega$, respectively.
\end{lemma}

\setcounter{equation}{0}
\setcounter{theorem}{0}

\section{Four tools}
\label{section3}

Let us first introduce a lemma concerning the unique determination of the inclusions with one Cauchy data for $N\ge2$ dimensions. We modify the proof which is given for the conductivity equation in \cite[Theorem 3.3, p. 72]{AmKa2007Springer}.  


\begin{lemma} 
\label{le: the unique determination} Let $\Omega$ be a bounded $C^2$ domain in $\mathbb R^N\ (N \ge 2)$ with boundary $\partial\Omega$, and let $D_1$ and $D_2$ be two bounded Lipschitz open sets, each of which may have finitely many connected components,  such that $D_1 \subset D_2 \subset \overline{D_2} \subset \Omega$ and both $\Omega\setminus\overline{D_1}$ and $\Omega\setminus\overline{D_2}$ are connected. $D_1$ and $D_2$ can be empty. Let $\sigma_j = \sigma_j(x)\ (j=1,2)$ be given by
$$
\sigma_j =
\begin{cases}
\sigma_c \quad&\mbox{in } D_j, \\
\sigma_s \quad&\mbox{in } \Omega \setminus D_j,
\end{cases}
$$
where $\sigma_c, \sigma_s$ are positive constants with $\sigma_c \not=\sigma_s$.  Let $g \in L^2(\partial\Omega)$ be a non-zero function, and let $v_j =v_j(x) \in H^1(\Omega)\ (j=1,2)$ satisfy
\begin{equation}
\label{modified poisson equation}
\mbox{\rm div}(\sigma_j\nabla v_j) = v_j -1\ \mbox{ in } \Omega\ \mbox{ and }\  \sigma_s \frac {\partial v_j}{\partial \nu} = g \ \mbox{ on } \partial\Omega,
\end{equation}
where $\nu$ denotes the unit outward normal vector to $\partial\Omega$. 
Then, if $v_1=v_2$  on $\partial\Omega$, $v_1=v_2$ in $\Omega$  and $D_1=D_2$.
\end{lemma} 


 \begin{remark}
 \label{poisson equation case}
 In {\rm \cite{Strieste2016}}, we dealt with Poisson's equation $\mbox{\rm div}(\sigma_j\nabla v_j) = -1$ instead of the first equation in \eqref{modified poisson equation}.
 If we replace the equation in \eqref{modified poisson equation} with Poisson's equation, then Lemma \ref{le: the unique determination} does not hold.
 Indeed,  denote by $B_s(x)$ the open ball  in $\mathbb R^N$ with radius $s > 0$ centered at $x \in \mathbb R^N$. For $x_0 \in \mathbb R^N$ and $\rho \in (0,1)$,  set 
 $$
 \Omega = B_1(x_0)\ \mbox{  and } D = B_\rho(x_0).
$$
Define two functions $u = u(r)$ and $v = v(r)$  by
\begin{eqnarray*}
&&u(r) = \frac 1{2N\sigma_s}(1-r^2)\ \mbox{ for } r \in [0,1],
\\
&& v(r) =  \left\{\begin{array}{rll}
 u(r) \ &\mbox{ for  }\  r \in [\rho,1],
\\
 \frac {\sigma_s}{\sigma_c}(u(r)-u(\rho)) + u(\rho)\ &\mbox{ for }\  r \in [0,\rho).
\end{array}\right.
\end{eqnarray*}
Then $v= v(|x-x_0|)$ satisfies 
$$
\mbox{\rm div}(\sigma\nabla v) = -1\ \mbox{ in } \Omega\   \mbox{ where }\ 
\sigma=
\begin{cases}
\sigma_c \quad&\mbox{in } D, \\
\sigma_s \quad&\mbox{in } \Omega \setminus D.
\end{cases}
$$
Since $\rho \in (0,1)$ can be chosen arbitrarily, the inclusion $D$ is not uniquely determined although  the solution $v$ is the same as $u$ outside $D$.
By the way this solution $v$ plays a key role in {\rm \cite{Strieste2016}}, but in the present paper we cannot use this function $v$ due to some technical reasons.
\end{remark}

\noindent
{\it Proof of Lemma \ref{le: the unique determination}.\ } For every $\eta \in H^1(\Omega)$ and for $j=1, 2$, we have
\begin{equation}
\label{weak form of Neumann problem}
\int_\Omega \left\{\sigma_j\nabla v_j\cdot\nabla \eta+ v_j\eta\right\} dx = \int_\Omega\eta\ dx + \int_{\partial\Omega}g\eta \ dS_x,
\end{equation}
where $dS_x$ denotes the area element. Hence it follows that for every $\eta \in H^1(\Omega)$
\begin{equation}
\label{weak difference}
\int_\Omega \left\{\sigma_1\nabla(v_1-v_2)\cdot\nabla \eta+ (v_1-v_2)\eta\right\} dx = (\sigma_c-\sigma_s)\int_{D_2\setminus D_1} \nabla v_2\cdot\nabla \eta\ dx.
\end{equation}
Substituting $\eta \equiv 1$ in \eqref{weak difference} yields that
\begin{equation}
\label{the same integral value}
\int_{\Omega}(v_1-v_2)\ dx = 0.
\end{equation}
Therefore, since $v_1-v_2 \in H^1_0(\Omega)$ and $\mbox{\rm div}(\sigma_1\nabla v_1) = v_1 -1\ \mbox{ in } \Omega$, we have from \eqref{the same integral value}
\begin{equation}
\label{for test function as distribution}
\int_\Omega \left\{\sigma_1\nabla v_1\cdot\nabla(v_1-v_2) + v_1(v_1-v_2)\right\} dx = \int_\Omega(v_1-v_2)\ dx  = 0
\end{equation}
By substituting $\eta = v_1$ in \eqref{weak difference}, we obtain from \eqref{for test function as distribution}
$$
(\sigma_c-\sigma_s)\int_{D_2\setminus D_1} \nabla v_2\cdot\nabla v_1\ dx = 0.
$$
Hence, substituting $\eta = v_1-v_2$ in \eqref{weak difference} yields that
\begin{equation}
\label{case where core has larger conductivity}
\int_\Omega\left\{\sigma_1|\nabla(v_1-v_2)|^2 + (v_1-v_2)^2\right\} dx +(\sigma_c-\sigma_s)\int_{D_2\setminus D_1} |\nabla v_2|^2\ dx = 0.
 \end{equation}
Thus, if $\sigma_c > \sigma_s$, then $v_1= v_2$ in $\Omega$. Suppose that $D_1 \not= D_2$. Since $D_1$ and $D_2$ are bounded Lipschitz open sets, then $D_2\setminus D_1$ contains an interior point. Hence $v_1=v_2 \equiv 1$ in $\overline{D_2}\setminus D_1$. Moreover, since $\sigma_s\Delta v_1 = v_1-1$ in $\Omega \setminus \overline{D_1}$, the function $v_1$ is real analytic in $\Omega \setminus \overline{D_1}$. Therefore, since $\Omega \setminus \overline{D_1}$ is connected,
  $v_1 \equiv 1$ in $\Omega \setminus \overline{D_1}$. This contradicts the assumption that $g\not=0$. Consequently, $D_1 = D_2$. If $\sigma_c < \sigma_s$, we interchange the roles of $v_1$ and $v_2$ to arrive at
$$
 \int_\Omega\left\{\sigma_2|\nabla(v_1-v_2)|^2 + (v_1-v_2)^2\right\} dx +(\sigma_s-\sigma_c)\int_{D_2\setminus D_1} |\nabla v_1|^2\ dx = 0,
$$
 which yields the same conclusion. \qed 
 

Also, we give a comparison lemma for partial differential inequalities.
 
 \begin{lemma} 
\label{le: the integral comparison} Let $\Omega$ be a bounded $C^2$ domain in $\mathbb R^N\ (N \ge 2)$ with boundary $\partial\Omega$, and let $D_1$ and $D_2$ be two bounded Lipschitz open sets, each of which may have finitely many connected components,  such that $D_1 \subset D_2 \subset \overline{D_2} \subset \Omega$ and both $\Omega\setminus\overline{D_1}$ and $\Omega\setminus\overline{D_2}$ are connected. $D_1$ and $D_2$ can be empty. Let $\sigma_j = \sigma_j(x)\ (j=1,2)$ be given by
$$
\sigma_j =
\begin{cases}
\sigma_c \quad&\mbox{in } D_j, \\
\sigma_s \quad&\mbox{in } \Omega \setminus D_j,
\end{cases}
$$
where $\sigma_c, \sigma_s$ are positive constants with $\sigma_c  \not= \sigma_s$.  Let $v_j =v_j(x) \in H^1(\Omega)\ (j=1,2)$ satisfy
$$
\mbox{\rm div}(\sigma_1\nabla v_1) \le v_1-1\ \mbox{ and }\  \mbox{\rm div}(\sigma_2\nabla v_2) \ge v_2-1\mbox{ in } \Omega.
$$
If $v_1\ge v_2$  on $\partial\Omega$ and $(\sigma_c\nabla v_2 - \sigma_s\nabla v_1)\cdot(\nabla v_2 - \nabla v_1) \ge 0$ in $D_2\setminus \overline{D_1}$,  then
$v_1 \ge v_2$ in $\Omega$.
\end{lemma} 

\noindent
{\it Proof.\ } Set $w = (v_2-v_1)^+(= \max\{v_2-v_1,0\})$. Since  $v_1\ge v_2$  on $\partial\Omega$, $w \in H^1_0(\Omega)$ and $w \ge 0$ in $\Omega$.
Therefore we have
$$
-\int_\Omega \sigma_1\nabla v_1\cdot\nabla w\ dx \le \int_\Omega(v_1w-w)\ dx\ \mbox{ and }\ -\int_\Omega \sigma_2\nabla v_2\cdot\nabla w\ dx \ge \int_\Omega(v_2w-w)\ dx.
$$
Thus
$$
\int_\Omega(\sigma_2\nabla v_2-\sigma_1\nabla v_1)\cdot \nabla w + \int_\Omega w^2 dx \le 0,
$$
and hence
$$
\int_{D_1}\!\!\sigma_c|\nabla w|^2 dx + \int_{\Omega\setminus D_2}\!\! \sigma_s|\nabla w|^2 dx + \int_\Omega w^2 dx + \int_{D_2 \setminus D_1}\!\! (\sigma_c\nabla v_2-\sigma_s\nabla v_1)\cdot \nabla w \ dx \le 0.
$$
This concludes that $w = 0$ in $\Omega$, since the fourth term in the left-hand side is nonnegative from the assumption.  \qed


Let us introduce a lemma for ordinary differential equations which can be proved by the power series  method and D'Alembert's method of reduction of order. 
\begin{lemma} 
\label{le: fundamental solutions of ODE}
For a number $\sigma > 0$ and $N \ge 2$, consider the ordinary differential equation for $f = f(r)$
\begin{equation}
\label{2nd order ODE for r}
f^{\prime\prime} + \frac {N-1}r f^\prime -\frac 1\sigma f = 0\ \mbox{ for } r \in (0,\infty).
\end{equation}
Then a fundamental set of solutions on the whole interval $(0,\infty)$ consists of two solutions $f_{reg}=f_{reg}(r)$ and $f_{sing}=f_{sing}(r)$ of the form:
\begin{equation}
\label{singular and regular solutions}
f_{sing}(r) = f_{reg}(r)\!\! \int_1^r\!\! \frac 1{s^{N-1}(f_{reg}(s))^2} ds\ \mbox{ and }\  f_{reg}(r) = \sum_{k=0}^\infty \frac {(N-2)!!}{[N+2(k-1)]!!\ k! \ 2^k \sigma^k}r^{2k}.
\end{equation}
Moreover, 
\begin{equation}
\label{limits at r equals zero}
\lim_{r \to +0} r^{N-1}f_{sing}^\prime(r) = 1,\ \lim_{r \to +0} f_{sing}(r) = -\infty,\ f_{reg}^\prime(0) = 0\mbox{ and }\  f_{reg}(0) = 1.
\end{equation}
Additionally, for every solution $f$ of \eqref{2nd order ODE for r} and every $\rho > 0$, the following formulas hold:
\begin{eqnarray}
&&f^\prime_{sing}(\rho)f_{reg}(\rho) -f_{sing}(\rho)f^\prime_{reg}(\rho) (= \rho^{1-N}) > 0; \label{sign of Wronskian}
\\
&&\mbox{ if }\ f = c_1f_{sing}+ c_2f_{reg}\ \mbox{ for two constants } c_1 \mbox{ and } c_2, \nonumber
\\
&&\mbox{ then }\ c_1 = \frac{f^\prime(\rho)f_{reg}(\rho) -f(\rho)f^\prime_{reg}(\rho)}{f^\prime_{sing}(\rho)f_{reg}(\rho) -f_{sing}(\rho)f^\prime_{reg}(\rho)}. \label{relationship between Wronskian and c_1}
 \end{eqnarray}
\end{lemma}

\noindent
{\it Proof.\ }  A simple application of the power series  method gives $f_{reg}$, and  D'Alembert's method of reduction of order gives $f_{sing}$. Thus \eqref{singular and regular solutions} holds, and hence both \eqref{limits at r equals zero} and   \eqref{sign of Wronskian} follow directly from \eqref{singular and regular solutions}.
    \eqref{sign of Wronskian} guarantees \eqref{relationship between Wronskian and c_1}. \qed


Finally, we give a comparison lemma for two solutions of \eqref{2nd order ODE for r} for different $\sigma$'s.
 \begin{lemma} 
\label{le: a comparison of ODE}
Let $0 < \sigma_1 < \sigma_2,\ N \ge 2$ and let $f_j = f_j(r)\ (j=1, 2)$ solve \eqref{2nd order ODE for r} with $\sigma = \sigma_j\ (j=1, 2)$, respectively.  Suppose that $f_1(\rho) = f_2(\rho)$ for some $\rho > 0$. Then the following assertions hold:
\begin{enumerate}[\rm (1)]
\item Assume that $\sigma_1f_1^\prime(\rho) = \sigma_2f_2^\prime(\rho) >0$. Then we have
\begin{itemize}
\item[{\rm (i)}] If there exists $s \in (0, \rho)$ such that  $f_1(s) = f_2(s)$ and $f_1(r) < f_2(r)$ for every $r \in (s,\rho)$, then $f_1^\prime(s) < 0$ and $f_2^\prime(s) < 0$.
 
 \item[{\rm (ii)}] If there exists $\ell \in (\rho, \infty)$ such that  $f_1(\ell) = f_2(\ell)$ and $f_1(r) > f_2(r)$ for every $r \in (\rho,\ell)$, then $f_1^\prime(\ell) < 0$ and $f_2^\prime(\ell) < 0$.
 \end{itemize}
\item Assume that $\sigma_1f_1^\prime(\rho) = \sigma_2f_2^\prime(\rho) <0$. Then we have

\begin{itemize}
\item[{\rm (i)}] If there exists $s \in (0, \rho)$ such that  $f_1(s) = f_2(s)$ and $f_1(r) > f_2(r)$ for every $r \in (s,\rho)$, then $f_1^\prime(s) > 0$ and $f_2^\prime(s) > 0$.
 
 \item[{\rm (ii)}] If there exists $\ell \in (\rho, \infty)$ such that  $f_1(\ell) = f_2(\ell)$ and $f_1(r) < f_2(r)$ for every $r \in (\rho,\ell)$, then $f_1^\prime(\ell) > 0$ and $f_2^\prime(\ell) > 0$.
 \end{itemize}
 \item Assume that $f_1^\prime(\rho) = f_2^\prime(\rho) =0$ and  $f_1(\rho) = f_2(\rho) < 0$. Then we have
\begin{itemize}
\item[{\rm (i)}] If there exists $s \in (0, \rho)$ such that  $f_1(s) = f_2(s)$ and $f_1(r) < f_2(r)$ for every $r \in (s,\rho)$, then $f_1^\prime(s) < 0$ and $f_2^\prime(s) < 0$.
 
 \item[{\rm (ii)}]  If there exists $\ell \in (\rho, \infty)$ such that  $f_1(\ell) = f_2(\ell)$ and $f_1(r) < f_2(r)$ for every $r \in (\rho,\ell)$, then $f_1^\prime(\ell) > 0$ and $f_2^\prime(\ell) > 0$.
 \end{itemize}
\end{enumerate}
\end{lemma}

\noindent
{\it Proof.\ } Observe that
 \begin{equation}
 \label{equation of the difference for 3}
\left(\sigma_1r^{N-1}f_1^\prime(r)\right)^\prime - \left(\sigma_2r^{N-1}f_2^\prime(r)\right)^\prime = r^{N-1}(f_1(r)-f_2(r))\ \mbox{ for } r >0.
\end{equation}
Let us first consider (3).   Since  $f_1^\prime(\rho) = f_2^\prime(\rho) =0,\ f_1(\rho) = f_2(\rho) < 0$ and  $0 < \sigma_1 < \sigma_2$, we observe that
$$
f_1^{\prime\prime}(\rho) = \frac 1{\sigma_1} f_1(\rho) < \frac 1{\sigma_2} f_2(\rho) = f_2^{\prime\prime}(\rho), 
$$
 and hence there exists a number $\delta > 0$ such that
 $$
 f_1(r) < f_2(r)\ \mbox{ for every } r \in (\rho-\delta, \rho) \cup (\rho, \rho + \delta).
 $$
 Let us prove (i).  
Since $f_1^\prime(\rho) = f_2^\prime(\rho) =0,\  f_1(s) = f_2(s)$ and $f_1(r) < f_2(r)$ for every $r \in (s,\rho)$,  
by integrating \eqref{equation of the difference for 3} over the interval $[s,\rho]$,  we have 
$$
f_1^\prime(s)\le f_2^\prime(s)\ \mbox{ and } - \left(\sigma_1s^{N-1}f_1^\prime(s) - \sigma_2s^{N-1}f_2^\prime(s)\right) = \int_s^\rho r^{N-1}(f_1(r)-f_2(r))\ dr < 0.
$$
These yield that $f_1^\prime(s) < 0$ and $f_2^\prime(s) < 0$,  since $0 < \sigma_1 < \sigma_2$. (ii) is proved similarly.

Let us consider (1).  Since $\sigma_1f_1^\prime(\rho) = \sigma_2f_2^\prime(\rho) > 0, \ f_1(\rho) = f_2(\rho)$ and $0 < \sigma_1 < \sigma_2$,  we observe that
$$
f_1^{\prime}(\rho) > f_2^{\prime}(\rho), 
$$
and hence there exists  a number $\delta > 0$ such that
 $$
 f_1(r) < f_2(r)\ \mbox{ for every } r \in (\rho-\delta, \rho)\ \mbox{ and  }\   f_1(r) >  f_2(r)\ \mbox{ for every } r \in(\rho, \rho + \delta).
 $$
Let us prove (i).  Since $\sigma_1f_1^\prime(\rho) = \sigma_2f_2^\prime(\rho),\  f_1(s) = f_2(s)$ and $f_1(r) < f_2(r)$ for every $r \in (s,\rho)$,  
by integrating \eqref{equation of the difference for 3} over the interval $[s,\rho]$,  we have 
$$
f_1^\prime(s)\le f_2^\prime(s)\ \mbox{ and } - \left(\sigma_1s^{N-1}f_1^\prime(s) - \sigma_2s^{N-1}f_2^\prime(s)\right) = \int_s^\rho r^{N-1}(f_1(r)-f_2(r))\ dr < 0.
$$
These yield that $f_1^\prime(s) < 0$ and $f_2^\prime(s) < 0$,  since $0 < \sigma_1 < \sigma_2$. (ii) is proved similarly.

Let us consider (2).  Since $\sigma_1f_1^\prime(\rho) = \sigma_2f_2^\prime(\rho) < 0, \ f_1(\rho) = f_2(\rho)$ and $0 < \sigma_1 < \sigma_2$,  we observe that
$$
f_1^{\prime}(\rho) < f_2^{\prime}(\rho), 
$$
and hence there exists  a number $\delta > 0$ such that
 $$
 f_1(r) > f_2(r)\ \mbox{ for every } r \in (\rho-\delta, \rho)\ \mbox{ and  }\   f_1(r) <  f_2(r)\ \mbox{ for every } r \in(\rho, \rho + \delta).
 $$
 Thus the conclusion follows from the same argument as in (1). \qed

\setcounter{equation}{0}
\setcounter{theorem}{0}

\section{Proof of Theorem \ref{th:stationary isothermic cauchy in two dimensions}}
\label{section4}

Let $u$ be the solution of  problem  \eqref{heat Cauchy} in two dimensions. For assertion (b) of Theorem \ref{th:stationary isothermic cauchy in two dimensions},
Lemma \ref{le: constant weingarten curvature} yields that $\gamma$ and $\Gamma$ are concentric circles. Denote by $x_0 \in \mathbb R^2$ the common center of $\gamma$ and $\Gamma$. By combining the initial  condition of problem  \eqref{heat Cauchy}  and  the assumption \eqref{stationary isothermic surface partially} with the real analyticity in $x$ of $u$ over $\mathbb R^2 \setminus\overline{D}$ coming from $\sigma_s=\sigma_m$,  we see that $u$ is radially symmetric with respect to $x_0$  in $x$ on $\left(\mathbb R^2\setminus\overline{D}\right) \times (0,\infty)$. Here we used the assumption that $\Omega\setminus\overline{D}$ is connected. Moreover, in view of the initial  condition of problem  \eqref{heat Cauchy},  we can distinguish the following two cases:
$$
\mbox{\rm (I)  } \Omega \mbox{ is a disk;}\qquad \mbox{\rm (II)  } \Omega \mbox{ is an annulus.}
$$ 
For assertion (a) of Theorem \ref{th:stationary isothermic cauchy in two dimensions}, Lemma \ref{le: constant weingarten curvature 2} yields that $\partial G$ and $\partial\Omega$ are concentric circles, since every component of $\partial\Omega$ is a circle with the same curvature. Therefore, only the case (I) remains for assertion (a) of Theorem \ref{th:stationary isothermic cauchy in two dimensions}. Also, denoting by $x_0 \in \mathbb R^2$ the common center of $\partial G$ and $\partial\Omega$ and combining the initial  condition of problem  \eqref{heat Cauchy}  and  the assumption \eqref{stationary isothermic surface} with the real analyticity in $x$ of $u$ over $\Omega \setminus\overline{D}$ yield that  $u$ is radially symmetric with respect to $x_0$  in $x$ on $\left(\mathbb R^2\setminus\overline{D}\right) \times (0,\infty)$. 

By virtue of (a) of Lemma \ref{le:initial behavior exponential decay and decay at infinity},  we can introduce the following three auxiliary functions $U = U(x), \ V = V(x)$ and $W =W(x)$ by
\begin{eqnarray}
&& U(x) = \int_0^\infty e^{-t}(1-u(x,t) )\ dt\quad \mbox{ for } x \in \Omega \setminus D, \label{auxiliary function on shell 2}
\\
&& V(x) = \int_0^\infty e^{-t}(1-u(x,t) )\ dt\quad \mbox{ for } x \in D,\label{auxiliary function on core 2}
\\
&& W(x) = \int_0^\infty e^{-t}(1-u(x,t) )\ dt\quad \mbox{ for } x \in \mathbb R^2 \setminus \Omega.\label{auxiliary function on medium}
\end{eqnarray}
Then we observe that
\begin{eqnarray}
&&0< U < 1\mbox{ in } \Omega \setminus\overline{D},\ 0< V < 1 \mbox{ in } D,\ 0 < W < 1\mbox{ in } \mathbb R^2\setminus\overline{\Omega}, \label{all bounded from above and below}
\\
&&\sigma_s\Delta U =U - 1  \mbox{ in } \Omega \setminus\overline{D},\ \sigma_c\Delta V =V- 1 \mbox{ in } D,\ 
  \sigma_m\Delta W =W \mbox{ in } \mathbb R^2\setminus\overline{\Omega}, \label{poisson and Laplace equations}
\\
&&U = V\ \mbox{ and }\  \sigma_s \frac {\partial U}{\partial\nu} = \sigma_c \frac {\partial V}{\partial\nu}  \ \mbox{ on } \partial D, \label{transmission condition between U and V}
\\
&&U = W\ \mbox{ and }\ \sigma_s \frac {\partial U}{\partial\nu} = \sigma_m \frac {\partial W}{\partial\nu}  \ \mbox{ on } \partial\Omega, \label{transmission condition between U and W}
\\
&& \lim_{|x| \to \infty} W(x) = 0, \label{decay at infinity}
\end{eqnarray}
where $\nu = \nu(x)$ denotes the unit outward normal vector to $\partial D$ at $x \in \partial D$ or to $\partial\Omega$ at $x \in \partial\Omega$, and the transmission conditions on $\partial D$ or on $\partial\Omega$  are given by \eqref{transmission condition between U and V}  or  by \eqref{transmission condition between U and W}, respectively. Here we used Lemma \ref{le:initial behavior exponential decay and decay at infinity} together with Lebesgue's dominated convergence theorem to obtain \eqref{decay at infinity}.

\vskip 2ex

 We first show that case (II)  for  assertion (b) of Theorem \ref{th:stationary isothermic cauchy in two dimensions} never occurs.  Set 
 $$
 \Omega = B_{\rho_+} \setminus \overline{B_{\rho_-}} \mbox{ with } B_{\rho_+} = B_{\rho_+}\!(x_0)  \mbox{ and } B_{\rho_-} = B_{\rho_-}\!(x_0) 
 $$
 for some numbers $\rho_+ > \rho_- > 0$.
 Since $u$ is radially symmetric with respect to $x_0$  in $x$ on $\left(\mathbb R^2\setminus\overline{D}\right) \times (0,\infty)$,  $W$ is radially symmetric with respect to $x_0$. 
 Observe from \eqref{all bounded from above and below} and \eqref{poisson and Laplace equations}  that 
 $$
 \Delta W > 0\   \mbox{ in }  B_{\rho_-} \cup \left(\mathbb R^2\setminus \overline{B_{\rho_+}}\right).
 $$
 Therefore,  in view of \eqref{decay at infinity},  by applying the strong maximum principle to the radially symmetric function $W$, we see that the positive maximum values $\max\limits_{\overline{B_{\rho_-}}} W$ and $\max\limits_{\mathbb R^2 \setminus B_{\rho_+}} W$ are achieved  only on $\partial B_{\rho_-}$ and on $\partial B_{\rho_+}$, respectively.  Hence,  Hopf's boundary point lemma yields that
 \begin{equation}
 \label{signs of radial derivatives}
 W^\prime(\rho_-) > 0\ \mbox{ and } W^\prime(\rho_+) < 0, 
 \end{equation}
 where we write $W^\prime = \frac d{dr}W$ for $r =|x-x_0|$. See \cite[Lemma 3.4, p. 34]{GT1983} for Hopf's boundary point lemma.

 Also, we see that  $U-1$  solves the ordinary differential equation \eqref{2nd order ODE for r} with $\sigma=\sigma_s$ in $r =|x-x_0|$. Moreover,  since $\Omega \setminus \overline{D}$ is connected, $U-1$ is extended as a solution of \eqref{2nd order ODE for r} for all $r =|x-x_0|$ in $\mathbb R^2 \setminus \{x_0\}$.  Thus,  it follows from \eqref{all bounded from above and below},  \eqref{poisson and Laplace equations}  and \eqref{transmission condition between U and W} 
that for $r =|x-x_0| \ge 0$
\begin{eqnarray}
&& \sigma_s\left(U^{\prime\prime} + \frac 1rU^\prime\right)= U-1 < 0 \  \mbox{ for } \rho_-\le r \le \rho_+, \label{some monotonicity of U}
\\
&& U^\prime(\rho_-) > 0\ \mbox{ and } U^\prime(\rho_+) < 0. \label{sign of derivatives of U on the boundary of Omega}
\end{eqnarray}
We set $D_1=\emptyset$ and $D_2 = D$, and we consider two functions $v_j= v_j(x)\in H^1(\Omega)\ (j=1,2)$ defined by 
\begin{equation}
\label{two functions in Lemma 3.1}
v_1 = U \ \mbox{ in } \Omega\ \mbox{ and } v_2 = \left\{\begin{array}{rll}
 U \ &\mbox{ in }\ \Omega \setminus D,
\\
 V \ &\mbox{ in }\  D.
\end{array}\right.
\end{equation}
Then we apply Lemma \ref{le: the unique determination} to these two functions $v_j= v_j(x)\in H^1(\Omega)$ to see that $v_1=v_2$ in $\Omega$  and $\emptyset=D$, which is a contradiction. 
Thus case (II)  for  assertion (b) of Theorem \ref{th:stationary isothermic cauchy in two dimensions} never occurs. 

\vskip 2ex
It remains to consider case (I). Set 
$$
\Omega = B_{\rho_0}(x_0)
$$
for some number $\rho_0 > 0$. Since $u$ is radially symmetric with respect to $x_0$  in $x$ on $\left(\mathbb R^2\setminus\overline{D}\right) \times (0,\infty)$,  $W$ is radially symmetric with respect to $x_0$. 
 Observe from \eqref{all bounded from above and below} and \eqref{poisson and Laplace equations}  that 
 $$
 \Delta W > 0\   \mbox{ in }  \mathbb R^2\setminus \overline{B_{\rho_0}(x_0)}.
 $$
 Therefore,  in view of \eqref{decay at infinity},  by applying the strong maximum principle to the radially symmetric function $W$, we see that the positive maximum value $\max\limits_{\mathbb R^2 \setminus B_{\rho_0}(x_0)} W$ is achieved  only on $\partial B_{\rho_0}(x_0)$.  Hence,  Hopf's boundary point lemma yields that
 \begin{equation}
 \label{sign of radial derivative}
 W^\prime(\rho_0) < 0.
 \end{equation}
 Since $U-1$  solves the ordinary differential equation \eqref{2nd order ODE for r} with $\sigma=\sigma_s$ in $r =|x-x_0|$ and $\Omega \setminus \overline{D}$ is connected, $U-1$ is extended as a solution of \eqref{2nd order ODE for r} for all $r =|x-x_0|$ in $\mathbb R^2 \setminus \{x_0\}$.  Therefore  it follows from Lemma \ref{le: fundamental solutions of ODE},  \eqref{all bounded from above and below}, \eqref{poisson and Laplace equations}, \eqref{transmission condition between U and W} and \eqref{sign of radial derivative}  that for $r =|x-x_0| \ge 0$
\begin{eqnarray}
&& U = c^*_1 f_{sing}(r) + c^*_2 f_{reg}(r) + 1 \  \mbox{ for } r > 0, \label{structure of U}
\\
&& \sigma_s\left(rU^\prime\right)^\prime = r (U-1) < 0 \  \mbox{ in  } \Omega \setminus \overline{D}, \label{some monotonicity of U 2}
\\
&&  U^\prime(\rho_0) < 0, \label{sign of derivative of U on the boundary of Omega 2}
\end{eqnarray}
where $c^*_1$ and $c^*_2$ are some constants and we chose $\sigma = \sigma_s$ in Lemma \ref{le: fundamental solutions of ODE}. 
We distinguish the following three cases:
$$
\mbox{\rm (i)  } c^*_1 = 0;\qquad \mbox{\rm (ii)  } c^*_1 < 0; \qquad \mbox{\rm (iii)  } c^*_1 > 0.
$$ 
Let us consider case (i) first. Notice that $U$ is smooth at $x=x_0$. Then, as in case (II), we set $D_1=\emptyset$ and $D_2 = D$, and define $v_j= v_j(x)\in H^1(\Omega)\ (j=1,2)$ by \eqref{two functions in Lemma 3.1}. Thus, by applying Lemma \ref{le: the unique determination} to these two functions $v_j= v_j(x)\in H^1(\Omega)$ to see that $v_1=v_2$ in $\Omega$  and $\emptyset=D$, which is a contradiction. Hence case (i) never occurs.

\vskip 2ex
Let us proceed to case (ii). Then it follows from Lemma \ref{le: fundamental solutions of ODE} that
\begin{equation}
\label{limits at zero of U}
\lim_{r \to 0}U(r) = -\lim_{r \to 0} U^\prime(r) =+\infty\ \mbox{ and } x_0\in D.
\end{equation}
Moreover, we notice that
\begin{equation}
\label{strictly decreasing of U}
U^\prime(r) < 0\ \mbox{ if } \rho_0 \ge r > 0.
\end{equation}
Indeed, by setting $h = U^\prime(r)$, we have
\begin{equation}
\label{the elliptic equation for the derivative}
\sigma_s\Delta h - \left(\frac {\sigma_s}{r^2} + 1\right) h = 0 \ \mbox{ in } B_{\rho_0}(x_0) \setminus\{ x_0 \}.
\end{equation}
Since $h$ is negative on $\partial B_{\rho_0}(x_0) \cup \partial B_\varepsilon(x_0)$ for sufficiently small $\varepsilon >0$,
the strong maximum principle yields that $h$ is negative in $B_{\rho_0}(x_0) \setminus \overline{B_\varepsilon(x_0)}$.

Let us choose the connected component $D_*$ of $D$ satisfying $x_0 \in D_*$. Then, since $\overline{D_*} \subset \Omega = B_{\rho_0}(x_0)$, we see that there exist $\rho_{*1},  \rho_{*2} \in  (0, \rho_0)$  and $x_{*1}, x_{*2} \in \partial D_*$ which  satisfy that $\rho_{*1} \le \rho_{*2}$ and
\begin{eqnarray}
&&U(\rho_{*1}) = \max\{ U(r) : r = |x-x_0|, x \in \partial D_*\} \mbox{ and } \rho_{*1} =|x_{*1}-x_0|, \label{maximum point}
\\
&&U(\rho_{*2}) = \min\{ U(r) : r = |x-x_0|, x \in \partial D_*\} \mbox{ and } \rho_{*2} =|x_{*2}-x_0|. \label{minimum point}
\end{eqnarray}
Notice that $\nu(x_{*i})$ equals $\frac{x_{*i}-x_0}{\rho_{*i}}$ for $i=1, 2$.  Also, the case where $\rho_{*1}=\rho_{*2}$ may occur for instance if $D_*$ is a disk centered at $x_0$.  When $\rho_{*1}=\rho_{*2}$, $D_*$ must be a disk centered at $x_0$ because of \eqref{strictly decreasing of U}. By setting $D_1=D_*$ and $D_2 = D$,  we consider two functions $v_j= v_j(x)\in H^1(\Omega)\ (j=1,2)$ defined by 
\begin{equation}
\label{two functions in Lemma 3.1-2}
v_1 = \left\{\begin{array}{rll}
 U \ &\mbox{ in }\ \Omega \setminus D_*,
\\
 V \ &\mbox{ in }\  D_*,
 \end{array}\right.
\ \mbox{ and } v_2 = \left\{\begin{array}{rll}
 U \ &\mbox{ in }\  \Omega \setminus D,
\\
 V \ &\mbox{ in }\  D.
\end{array}\right.
\end{equation}
Then we apply Lemma \ref{le: the unique determination} to these two functions $v_j= v_j(x)\in H^1(\Omega)$ to see that $v_1=v_2$ in $\Omega$  and $D_*=D$, which gives the desired conclusion of Theorem \ref{th:stationary isothermic cauchy in two dimensions}.  Hereafter in case (ii), we may assume that $\rho_{*1} < \rho_{*2}$.

Let $g_j = g_j(r) \ (j=1,2)$ be the unique solution of the Cauchy problem:
\begin{equation}
\label{two auxiliary radial solutions}
\sigma_c\left(g_j^{\prime\prime} + \frac 1r g_j^\prime\right) = g_j- 1\ \mbox{ for } r > 0,\ g_j(\rho_{*j}) = U(\rho_{*j})\ \mbox{ and } g_j^\prime(\rho_{*j}) =\frac {\sigma_s}{\sigma_c} U^\prime(\rho_{*j}).
\end{equation}
Then we observe that for $j=1, 2$
\begin{eqnarray}
&&\sigma_c\Delta V =V- 1 \ \mbox{ in } D_*, \ \sigma_c\Delta g_j =g_j- 1 \ \mbox{ in } D_* \setminus\{x_0\}, \label{PDEs in the components}
\\
&&g_j =V\ \mbox{ and } \frac {\partial g_j}{\partial \nu} = \frac {\partial V}{\partial \nu} \ \mbox{ at } x_{*j} \in \partial D_*. \label{boundary critical point j}
\end{eqnarray}

Let us distinguish the following two cases provided that $\rho_{*1} < \rho_{*2}$:
$$
\mbox{\rm (ii-a)  } \sigma_c < \sigma_s; \qquad \mbox{\rm (ii-b)  } \sigma_c > \sigma_s.
$$ 
In case (ii-a), we employ $g_2$. Since both $g_2-1$ and $U-1$ satisfy the ordinary differential equation \eqref{2nd order ODE for r} with $\sigma=\sigma_c$ and $\sigma=\sigma_s$ respectively and $g_2-1 = U-1 < 0$ at $r =\rho_{*2}$,  by taking  \eqref{strictly decreasing of U} and \eqref{limits at zero of U}
 into account we apply (2)-(i) of  Lemma \ref{le: a comparison of ODE} to these two solutions and conclude that
 \begin{equation}
\label{g against U}
g_2 \ge U (= V)  \mbox{ on } \partial D_*\ \mbox{ and } \lim_{x \to x_0}g_2= \infty.
\end{equation}
Thus it follows from \eqref{PDEs in the components}, \eqref{g against U} and the strong comparison principle  that
$$
g_2 > V \ \mbox{ in } D_* \setminus\{x_0\}.
$$
This contradicts \eqref{boundary critical point j} by Hopf's boundary point lemma. Thus case (ii-a) never occurs. 

\vskip 2ex
Let us proceed to case (ii-b).  We employ both $g_1$ and $g_2$. Since both $g_j-1$ and $U-1$ satisfy the ordinary differential equation \eqref{2nd order ODE for r} with $\sigma=\sigma_c$ and $\sigma=\sigma_s$ respectively and $g_j-1 = U-1 < 0$ at $r =\rho_{*j}$,  by taking  \eqref{strictly decreasing of U} and \eqref{limits at zero of U}
 into account we apply (2) of  Lemma \ref{le: a comparison of ODE} to these two solutions and conclude that the graphs of $g_j$ and $U$ intersect only at $r =\rho_{*j}$ in $(0, \rho_0)$
 for each $j =1, 2$.  By Lemma \ref{le: fundamental solutions of ODE}, we may set for each $j=1,2$
 \begin{equation}
 \label{representation of g_j by the fundamental solutions}
g_j(r) = c_{j, 1} f_{sing}(r) + c_{j, 2} f_{reg}(r) + 1 \  \mbox{ for } r > 0, 
\end{equation}
where  $c_{j, 1}$ and $c_{j, 2}$ are some constants and we chose $\sigma = \sigma_c$ in Lemma \ref{le: fundamental solutions of ODE}. When either $c_{1,1}$ or $c_{2,1}$ equals zero,  with the aid of Lemma \ref{le: the unique determination} we obtain a contradiction by employing the  argument similar to the case where $\rho_{*1}=\rho_{*2}$.
For instance, let us assume that $c_{2,1}$ equals zero. Then, by setting $D_1=D_*$ and $D_2 = B_{\rho_{*2}}(x_0)$,  we consider two functions $v_j= v_j(x)\in H^1(\Omega)\ (j=1,2)$ defined by 
\begin{equation}
\label{two functions in Lemma 3.1-3}
v_1= \left\{\begin{array}{rll}
 U\ &\mbox{ in }\  \Omega \setminus D_*,
\\
 V\ &\mbox{ in }\  D_*,
 \end{array}\right.
\ \mbox{ and } v_2= \left\{\begin{array}{rll}
 U \ &\mbox{ in }\  \Omega \setminus B_{\rho_{*2}}(x_0),
\\
 g_2\ &\mbox{ in }\ B_{\rho_{*2}}(x_0),
\end{array}\right.
\end{equation}
where $g_2 = g_2(|x-x_0|)$ for $x \in B_{\rho_{*2}}(x_0)$.
Therefore we apply Lemma \ref{le: the unique determination} to these two functions $v_j= v_j(x)\in H^1(\Omega)$ to see that $v_1=v_2$ in $\Omega$  and $D_*=B_{\rho_{*2}}(x_0)$, which contradicts the assumption that $\rho_{*1} < \rho_{*2}$. Thus we distinguish the following four cases:
$$
\mbox{\rm (ii-b-1)  } c_{1,1} >0, \ c_{2,1}>0; \qquad \mbox{\rm (ii-b-2)  } c_{1,1} <0, \ c_{2,1}<0;
$$
$$
\mbox{\rm (ii-b-3)  } c_{1,1} <0,\  c_{2,1}>0; \qquad \mbox{\rm (ii-b-4)  } c_{1,1} >0, \ c_{2,1}<0.
$$ 
The first three cases (ii-b-1), (ii-b-2), and (ii-b-3) never occur because of Hopf's boundary point lemma as in case (ii-a). For instance, in case (ii-b-1),  we employ $g_2$ and observe that 
\begin{equation}
\label{g_2 against U}
 g_2 \le U (= V)  \mbox{ on } \partial D_*\ \mbox{ and } \lim_{x \to x_0}g_2=  -\infty.
\end{equation}
Thus it follows from \eqref{PDEs in the components}, \eqref{g_2 against U} and the strong comparison principle  that
$$
g_2 < V \ \mbox{ in } D_* \setminus\{x_0\}.
$$
This contradicts \eqref{boundary critical point j} by Hopf's boundary point lemma. Thus case (ii-b-1) never occurs.  In case  (ii-b-2) we employ $g_1$, and in case (ii-b-3) we can employ 
either of $g_1$ and $g_2$.  Let us proceed to case (ii-b-4). In case (ii-b-4) we cannot employ either of them. 
Thus, for every $\rho \in (0,\rho_0)$, we consider  the unique solution 
$g_\rho=g_\rho(r)$ of the Cauchy problem:
\begin{equation}
\label{auxiliary radial solutions with continuous  parameter rho}
\sigma_c\left(g_\rho^{\prime\prime} + \frac 1r g_\rho^\prime\right) = g_\rho- 1\ \mbox{ for } r >0,\ g_\rho(\rho) = U(\rho)\  \mbox{ and } g_\rho^\prime(\rho) =\frac {\sigma_s}{\sigma_c} U^\prime(\rho).
\end{equation}
Note that $g_{\rho_{*j}} = g_j$ where $g_j\ (j=1, 2)$ are defined by \eqref{two auxiliary radial solutions}.
By Lemma \ref{le: fundamental solutions of ODE}, we may set for each $\rho \in (0,\rho_0)$
\begin{equation}
\label{representation of gs by the fundamental solutions}
g_\rho(r) = c_{1}(\rho) f_{sing}(r) + c_{2}(\rho) f_{reg}(r) + 1 \  \mbox{ for } r > 0, 
\end{equation}
where $c_{1}(\rho)$ and $c_{2}(\rho)$ are some constants and we chose $\sigma = \sigma_c$ in Lemma \ref{le: fundamental solutions of ODE}.
Note that  $c_{i}(\rho_{*j}) =  c_{j,i}$. It follows from formula \eqref{relationship between Wronskian and c_1} of Lemma \ref{le: fundamental solutions of ODE} and the definition of $g_\rho$ that
\begin{equation}
\label{represent c1 by derivatives of U}
c_{1}(\rho) = \frac{\frac {\sigma_s}{\sigma_c} U^\prime(\rho)f_{reg}(\rho) -(U(\rho)-1)f^\prime_{reg}(\rho)}{f^\prime_{sing}(\rho)f_{reg}(\rho) -f_{sing}(\rho)f^\prime_{reg}(\rho)}\ \mbox{ for each } \rho \in (0,\rho_0).
\end{equation}
In view of \eqref{limits at zero of U} and Lemma \ref{le: fundamental solutions of ODE}, we observe that there exists $\delta \in (0,\rho_{*1})$ satisfying 
\begin{equation}
\label{the situation to apply the intermediate value theorem}
c_1(\rho_{*2}) = c_{2,1}<0,\ c_1(\rho_{*1}) = c_{1,1} > 0,\ \mbox{ and }\  c_1(\rho) < 0\ \mbox{ if } \rho \in (0,\delta].
\end{equation}
Since  $c_{1}(\rho)$ is continuous in $\rho \in (0, \rho_0)$ because of \eqref{represent c1 by derivatives of U}, the intermediate value theorem yields that there exist two numbers $\rho_3$ and $\rho_4$ satisfying
$$
\delta < \rho_3 <  \rho_{*1} < \rho_4 < \rho_{*2}\ \mbox{ and }\ c_1(\rho_3) = c_1(\rho_4) = 0.
$$
Then, by setting $D_1=B_{\rho_3}(x_0)$ and $D_2 = B_{\rho_4}(x_0)$,  we consider two functions $v_j= v_j(x)\in H^1(\Omega)\ (j=1,2)$ defined by 
\begin{equation}
\label{two functions in Lemma 3.1-4}
v_1= \left\{\begin{array}{rll}
 U\ &\mbox{ in }\  \Omega \setminus B_{\rho_3}(x_0),
\\
 g_{\rho_3}\ &\mbox{ in }\  B_{\rho_3}(x_0),
 \end{array}\right.
\ \mbox{ and }\ v_2= \left\{\begin{array}{rll}
 U \ &\mbox{ in }\  \Omega \setminus B_{\rho_4}(x_0),
\\
 g_{\rho_4}\ &\mbox{ in }\ B_{\rho_4}(x_0),
\end{array}\right.
\end{equation}
where $g_{\rho_3} = g_{\rho_3}(|x-x_0|)$ for $x \in B_{\rho_3}(x_0)$ and $g_{\rho_4} = g_{\rho_4}(|x-x_0|)$ for $x \in B_{\rho_4}(x_0)$.
Therefore we apply Lemma \ref{le: the unique determination} to these two functions $v_j= v_j(x)\in H^1(\Omega)$ to see that $v_1=v_2$ in $\Omega$ and $B_{\rho_3}(x_0)=B_{\rho_4}(x_0)$, which is a contradiction. Thus case (ii-b-4) never occurs.

\vskip 2ex
It remains to consider case (iii) where $c^*_1 > 0$. Then it follows from Lemma \ref{le: fundamental solutions of ODE} that
\begin{equation}
\label{limits at zero of U 2}
\lim_{r \to 0}U(r) = -\lim_{r \to 0} U^\prime(r) =-\infty\ \mbox{ and } x_0\in D.
\end{equation}
Moreover, we notice that there exists a number $\rho_{max} \in (0,\rho_0)$ satisfying
\begin{equation}
\label{shape of U on the interval 2}
U^\prime > 0\ \mbox{ on } (0,\rho_{max}),\ U^\prime(\rho_{max}) = 0\ \mbox{ and }\  U^\prime < 0\ \mbox{ on } (\rho_{max},\rho_0].
\end{equation}
Indeed, because of \eqref{limits at zero of U 2} and \eqref{sign of derivative of U on the boundary of Omega 2} there exists at least one $\hat{\rho} \in (0,\rho_0)$ with $U^\prime(\hat{\rho}) = 0$. Hence, by setting $h = U^\prime(r)$, we have
\eqref{the elliptic equation for the derivative},  and hence for sufficiently small $\varepsilon > 0$ we apply the strong maximum principle to $h$ on $B_{\hat{\rho}}(x_0)\setminus \overline{B_\varepsilon(x_0)}$ and $B_{\rho_0}(x_0)\setminus \overline{B_{\hat{\rho}}(x_0)}$, respectively, to obtain \eqref{shape of U on the interval 2} with $\hat{\rho} = \rho_{max}$.
Here we eventually know that such a number $\hat{\rho}$ is unique and therefore we set $\hat{\rho} = \rho_{max}$, since $U$ achieves its maximum at $r = \rho_{max}$.

Let us choose the connected component $D_*$ of $D$ satisfying $x_0 \in D_*$.  Then, since $\overline{D_*} \subset \Omega = B_{\rho_0}(x_0)$,  as in case (ii), we see that there exist $\rho_{*1},  \rho_{*2} \in  (0, \rho_0)$  and $x_{*1}, x_{*2} \in \partial D_*$ which  satisfy \eqref{maximum point} and \eqref{minimum point}. In view of the shape of the graph of $U$, we have from the transmission condition \eqref{transmission condition between U and V} that at $x_{*i} \in \partial D_*, i=1, 2, $
\begin{equation}
\label{vanish normal derivative}
\frac {\partial V}{\partial\nu} = \frac{\sigma_s}{\sigma_c} \frac {\partial U}{\partial\nu}
 = \left\{\begin{array}{rll}
 0 \ &\mbox{ if } \rho_{*i} = \rho_{max}\ ,
 \\
\frac{\sigma_s}{\sigma_c} U^\prime \ &\mbox{ if } \rho_{*i} \not= \rho_{max}\ ,
\end{array}\right.
\end{equation}
where, in order to see that  $\nu(x_{*i})$  equals $\frac{x_{*i}-x_0}{\rho_{*i}}$ if $\rho_{*i} \not= \rho_{max} $, we used the fact that both $D_*$ and $B_{\rho_0}(x_0) \setminus \overline{D_*}$ are connected and $x_0 \in D_*$. Then we observe that for $j=1, 2$ both \eqref{PDEs in the components} and \eqref{boundary critical point j} hold also in case (iii).
Also, the case where $\rho_{*1}=\rho_{*2}$ may occur for instance if $D_*$ is a disk centered at $x_0$.   When $\rho_{*1}=\rho_{*2}$,  $D_*$ must be a disk centered at $x_0$ because of \eqref{shape of U on the interval 2}.  Hence, by employing the same argument as in the case where $\rho_{*1}=\rho_{*2}$ in case (ii) (see \eqref{two functions in Lemma 3.1-2}), Lemma \ref{le: the unique determination} yields $D =D_*$, which is the desired conclusion of Theorem \ref{th:stationary isothermic cauchy in two dimensions}.
Therefore,  hereafter we may assume that $\rho_{*1}\not=\rho_{*2}$. Then we notice that $\rho_{*2} > \rho_{max}$. Indeed, if $\rho_{*2} \le \rho_{max}$,  then $\rho_{*2} < \rho_{max}$ since $\rho_{*1}\not=\rho_{*2}$.  By \eqref{shape of U on the interval 2} and \eqref{vanish normal derivative}, $\frac {\partial V}{\partial\nu} > 0$ at $x_{*2} \in \partial D_*$. This implies that $V$ achieves its minimum over  $\overline{D_*}$ at some interior point in $D_*$, which contradicts the fact that $\Delta V < 0$ in $D_*$ because of the strong maximum principle. Since $\Omega\setminus \overline{D}$ is connected, $\partial D_*$ must be connected.  (Here $\partial D_*$ must be a simple closed curve in the plane. )  Distinguish two cases provided that $\rho_{*1}\not=\rho_{*2}$:
$$
\mbox{\rm (iii-1) }\  \partial D_* \cap \partial B_{\rho_{max}}\!(x_0) =\emptyset ; \qquad  \mbox{\rm (iii-2) }\  \partial D_* \cap \partial B_{\rho_{max}}\!(x_0) \not=\emptyset.
$$
 In case (iii-1) $\partial D_* \subset B_{\rho_0}(x_0) \setminus \overline{B_{\rho_{max}}\!(x_0) }$,  since $x_{*2} \in  \partial D_* \setminus \overline{B_{\rho_{max}}\!(x_0) }$ and $\partial D_*$ is connected.
In case (iii-2) $\rho_{*1}=\rho_{max}$ and $x_{*1} \in \partial D_* \cap \partial B_{\rho_{max}}\!(x_0) $.

Let us consider case (iii-1). We have that $\rho_{max} < \rho_{*1} < \rho_{*2} < \rho_0$ because of \eqref{shape of U on the interval 2}. Distinguish the following two cases:
$$
\mbox{\rm (iii-1-a)  } \sigma_c < \sigma_s; \qquad \mbox{\rm (iii-1-b)  } \sigma_c > \sigma_s.
$$ 
In case (iii-1-b), by employing $g_2$ and using the same argument as in case (ii-a) to obtain a contradiction by Hopf's boundary point lemma, we can see that case (iii-1-b) never occurs.
Here we used (2)-(ii) of Lemma \ref{le: a comparison of ODE} to obtain that $\lim\limits_{x \to x_0}g_2 = -\infty$. In case (iii-1-a), by employing all the functions $g_1, g_2$ and $g_\rho$ for $\rho_{\max} \le \rho <\rho_0$ and using the same argument as in case (ii-b) to obtain a contradiction, we can see that case (iii-1-a) never occurs. Here, instead of using that $c_1(\rho) < 0$ for $\rho \in (0,\delta]$ in \eqref{the situation to apply the intermediate value theorem} in case (ii-b), we used the fact that $c_1(\rho_{max})  > 0$.

Let us proceed to case (iii-2). Distinguish the following two cases:
$$
\mbox{\rm (iii-2-a)  } \sigma_c < \sigma_s; \qquad \mbox{\rm (iii-2-b)  } \sigma_c > \sigma_s.
$$ 
In case (iii-2-a), we employ $g_1$ and we have from (3) of Lemma \ref{le: a comparison of ODE} that  
\begin{equation}
\label{g_1 against U}
 g_1 \le U (= V)  \mbox{ on } \partial D_*\ \mbox{ and } \lim_{x \to x_0}g_1= -\infty.
\end{equation}
Then it follows from \eqref{PDEs in the components}, \eqref{g_1 against U} and the strong comparison principle that
$$
g_1< V\ \mbox{ in } D_* \setminus\{x_0\}.
$$
This contradicts \eqref{boundary critical point j} by Hopf's boundary point lemma. Thus case (iii-2-a) never occurs.

In case (iii-2-b), we employ $g_2$ and we have from (2) of Lemma \ref{le: a comparison of ODE} that
\begin{equation}
\label{g_2 against U at the origin}
\lim_{x \to x_0}g_2=  -\infty,
\end{equation}
and moreover it follows that if $\rho_{max} \le r < \rho_{*2}$ then
$$
U(r) > g_2(r)\ \mbox{ and }\ -r(\sigma_sU^\prime(r)-\sigma_cg_2^\prime(r)) = \int_r^{\rho_{*2}} t(U(t)-g_2(t)) dt > 0,
$$
which implies that
$$
\sigma_cg_2^\prime(r) > \sigma_s U^\prime(r).
$$
Also, since $U^\prime(r) \le 0$ and $\sigma_c > \sigma_s$, if $\rho_{max} \le r < \rho_{*2}$ then
$$
g_2^\prime(r) > U^\prime(r).
$$
Thus we have
\begin{equation}
\label{key inequality part 1}
g_2^\prime(\rho_{max}) > 0\ \mbox{ and }\ (\sigma_cg_2^\prime(r) - \sigma_s U^\prime(r))(g_2^\prime(r) - U^\prime(r)) > 0\ \mbox{ if } \rho_{max} \le r < \rho_{*2}.
\end{equation}
In view of Lemma \ref{le: fundamental solutions of ODE}, we can find a constant $\beta > 0$ to get
\begin{equation}
\label{the continuity at max}
1 -\beta f_{reg}(\rho_{max}) = g_2(\rho_{max}),
\end{equation}
where we chose $\sigma = \sigma_c$ in Lemma \ref{le: fundamental solutions of ODE}. Then we introduce a function $v_2 = v_2(r)$ given by
$$
v_2(r) = \left\{\begin{array}{rlll}
1 -\beta f_{reg}(r) \ &\mbox{ if }\ 0 \le r < \rho_{max},
\\
g_2(r)\ &\mbox{ if }\  \rho_{max} \le r < \rho_{*2},
\\
U(r)\ &\mbox{ if }\  \rho_{*2} \le r \le \rho_0.
 \end{array}\right.
 $$
 Hence we have in particular
 \begin{equation}
\label{key inequality part 2}
 (\sigma_cv_2^\prime(r) - \sigma_s U^\prime(r))(v_2^\prime(r) - U^\prime(r)) > 0\ \mbox{ if } 0< r < \rho_{max}.
\end{equation}
 Since $g_2^\prime(\rho_{max}) > 0$ and $f^\prime_{reg}(\rho_{max}) > 0$, with the aid of \eqref{the continuity at max} we know that
 $$
 \mbox{ div}(\sigma_2\nabla v_2) \ge v_2 -1\ \mbox{ in } \Omega,
 $$
 where we set $v_2 = v_2(|x-x_0|)$ for $x \in \Omega$ and
 $$
\sigma_2 =
\begin{cases}
\sigma_c \quad&\mbox{in } B_{\rho_{*2}}(x_0), \\
\sigma_s \quad&\mbox{in } \Omega \setminus B_{\rho_{*2}}(x_0).
\end{cases}
$$
 Moreover let us introduce a function $v_1 = v_1(x)$ given by
$$
v_1 = \left\{\begin{array}{rll}
V \ &\mbox{ in }\ D_*,
\\
U\ &\mbox{ in }\  \Omega \setminus D_*.
 \end{array}\right.
 $$
 Then 
 $$
 \mbox{ div}(\sigma_1\nabla v_1) = v_1 -1\ \mbox{ in } \Omega,
 $$
where we set
$$
\sigma_1 =
\begin{cases}
\sigma_c \quad&\mbox{in } D_*, \\
\sigma_s \quad&\mbox{in } \Omega \setminus D_*.
\end{cases}
$$
Therefore, since $D_* \subset B_{\rho_{*2}}(x_0)$, in view of \eqref{key inequality part 1} and \eqref{key inequality part 1} we can apply Lemma \ref{le: the integral comparison} to two open sets $D_1 =D_*$ and $D_2=B_{\rho_{*2}}(x_0)$ and we conclude that $v_1 \ge v_2$ in $\Omega$. Hence it follows from the strong comparison principle that in particular
$$
g_2 < V \ \mbox{ in } D_* \setminus\{x_0\}.
$$
This contradicts \eqref{boundary critical point j} by Hopf's boundary point lemma. Thus case (iii-2-b) never occurs.
The proof of Theorem \ref{th:stationary isothermic cauchy in two dimensions} is completed.

\setcounter{equation}{0}
\setcounter{theorem}{0}

\section{Concluding remarks and related two-phase elliptic overdetermined problems}
\label{section5}

As is mentioned in the end of section \ref{introduction}, the method employed in the present paper works also in $N (\ge 3)$ dimensions with the aid of the four key tools given in section \ref{section3}  and the preliminaries given in \cite[Section 2, pp. 169--180]{Strieste2016}, which are similar to those in section \ref{section2}. Hence  the same method as in the present paper also gives other proofs of Theorems \ref{th:stationary isothermic} and \ref{th:stationary isothermic cauchy} without using the explicit radially symmetric solutions of Poisson's equation over balls.  Moreover we can prove the following two theorems below concerning their related two-phase elliptic overdetermined problems.

To be precise, let $R > 0$ and consider the ball $B_R(0)$ in $\mathbb R^N (N \ge 2)$  with radius $R$ centered at the origin. Let $D$  be a bounded  $C^2$ open set in $\mathbb R^N$ which may have finitely many connected components, and assume that $B_R(0)\setminus\overline{D}$ is connected and $\overline{D} \subset B_R(0)$. Denote by $\sigma=\sigma(x)\ (x \in B_R(0))$  the conductivity distribution given by
$$
\sigma =
\begin{cases}
\sigma_c \quad&\mbox{in } D, \\
\sigma_s \quad&\mbox{in } B_R(0) \setminus D,
\end{cases}
$$
where $\sigma_c, \sigma_s$ are positive constants and $\sigma_c \not=\sigma_s$. Consider the unique solution $u \in H^1(B_R(0))$ of the following boundary value problem:
\begin{equation}
\label{modified poisson and Hermholtz equation}
\mbox{\rm div}(\sigma\nabla u) = \alpha u -\beta < 0\ \mbox{ in } B_R(0)\ \mbox{ and }\  u = c \ \mbox{ on } \partial B_R(0),
\end{equation}
where $\alpha \ge 0, \beta > 0$ and $c$ are real constants. Then we have the following two theorems:

\begin{theorem}
\label{th:constant Neumann boundary condition} Let $u$ be the solution of problem \eqref{modified poisson and Hermholtz equation}. Suppose that $u$ satisfies
\begin{equation}
\label{overdetermined Neumann 1}
\sigma_s\frac {\partial u}{\partial \nu} = d\ \mbox{ on } \partial B_R(0),
\end{equation}
where  $d$ is a negative constant and $\nu$ denotes the unit outward normal vector to $\partial B_R(0)$.  Then $D$ must be a ball centered at the origin.
\end{theorem}

\noindent
{\it Proof.\ } With the aid of the real analyticity of  the solution $u$ over $B_R(0)\setminus\overline{D}$,  assumption \eqref{overdetermined Neumann 1},  together with the uniqueness of the solution of the Cauchy problem for elliptic equations, yields that $u$ must be radially symmetric with respect to the origin over $B_R(0)\setminus\overline{D}$. Distinguish two cases:
$$
\mbox{\rm (i) } \alpha=0; \quad \mbox{\rm (ii) } \alpha > 0.
$$

In case (i), if we set
$$
\tilde{\sigma} = \frac \sigma \beta,
$$
then $u$ satisfies
$$
\mbox{\rm div}(\tilde{\sigma}\nabla u) = -1 < 0\ \mbox{ in } B_R(0)\ \mbox{ and } u^\prime(R) < 0.
$$
Hence, in order to conclude that $D$ is a ball centered at the origin, we can follow the proofs in \cite{Strieste2016} by using the explicit radially symmetric solutions of Poisson's equation over balls.  In fact, essentially this case has been proved in \cite{Strieste2016}, although the result is not stated in \cite{Strieste2016}.

In case (ii), if we set 
$$
\tilde{\sigma} = \frac \sigma \alpha\ \mbox{ and }\ v = \frac \alpha\beta u,
$$
then $v$ satisfies
$$
\mbox{\rm div}(\tilde{\sigma}\nabla v) = v-1 < 0\ \mbox{ in } B_R(0)\ \mbox{ and } v^\prime(R) < 0.
$$
Hence, in order to conclude that $D$ is a ball centered at the origin, we can follow the proof in section \ref{section4} of the present paper in $N (\ge 2)$ dimensions with the aid of the four key tools given in section \ref{section3}. \qed

\begin{theorem}
\label{th:constant Dirichlet boundary condition on the inner sphere} Let $u$ be the solution of problem \eqref{modified poisson and Hermholtz equation}. Suppose that  there exists $r \in (0,R)$ with $\overline{D} \subset B_r(0)$ and $u$ satisfies
\begin{equation}
\label{overdetermined Dirichlet 2}
u = d\ \mbox{ on } \partial B_r(0),
\end{equation}
where  $d$ is a constant with $ d > c$. Then $D$ must be a ball centered at the origin.
\end{theorem}

\noindent
{\it Proof.\ } By applying the maximum principle to the function $x_j\frac {\partial u}{\partial x_i}-x_i\frac{\partial u}{\partial x_j}$ for $i \not= j$ with assumption \eqref{overdetermined Dirichlet 2} we see that $u$ must be radially symmetric with respect to the origin over $B_R(0)\setminus\overline{B_r(0)}$ and hence the real analyticity of  the solution $u$ over $B_R(0)\setminus\overline{D}$ yields that $u$ must be radially symmetric with respect to the origin over $B_R(0)\setminus\overline{D}$. Moreover, it follows from the strong maximum principle and Hopf's boundary point lemma that
$u^\prime(R) < 0$. Then the rest of the proof runs along the same line as in that of Theorem \ref{th:constant Neumann boundary condition}. \qed

\end{document}